\documentclass[12pt]{amsart}
\usepackage[colorlinks=true, pdfstartview=FitV, linkcolor=blue, citecolor=blue, urlcolor=blue]{hyperref}
\usepackage{amssymb}
\calclayout

\def\C{\mathbb {C}}
\def\NN{{\mathcal N}}
\def\N{\mathbb {N}}
\def\Z{\mathbb {Z}}
\def\O{\mathcal O}
\def\lie#1{\mathfrak{ #1}}
\def\lieh{\lie h}
\def\lieg{\lie g}

\def\liet{\lie t}
\def\inv{^{-1}}
\newcommand{\Aut}{\operatorname{Aut}}

\newcommand{\SL}{\operatorname{SL}}
\newcommand{\Sp}{\operatorname{Sp}}
\newcommand{\SO}{\operatorname{SO}}
\newcommand{\Orth}{\operatorname{O}}
\newcommand{\PGL}{\operatorname{PGL}}

\def\GL{\operatorname{GL}}
\def\Hom{\operatorname{Hom}}
\def\End{\operatorname{End}}

\def\ad{\operatorname{ad}}

\def\phi{\varphi}

\def\quot#1#2{#1/\!\! /#2}
\def\I{{\mathcal I}}
\def\gr{\operatorname{gr}}
\def\sm#1#2#3#4{\left(\smallmatrix #1 & #2 \\ #3 & #4\endsmallmatrix\right)}

\numberwithin{equation}{subsection}

\newtheorem{theorem}[subsection]{Theorem}
\newtheorem{lemma}[subsection]{Lemma}
\newtheorem{proposition}[subsection]{Proposition}
\newtheorem{corollary}[subsection]{Corollary}
\newtheorem{conjecture}[subsection]{Conjecture}

\theoremstyle{definition}

\theoremstyle{remark}

\newtheorem{remark}[subsection]{Remark}

\newtheorem{example}[subsection]{Example}

\title[Linear maps preserving fibers]{\boldmath Linear maps preserving fibers} 
 \author{Gerald W. Schwarz}
\thanks{Partially supported by NSA Grant H98230-06-1-0023}
\address{Department of Mathematics\\
Brandeis University\\
Waltham, MA 02454-9110}
\email{schwarz@brandeis.edu}

\begin{document}
\begin{abstract}
Let $G\subset\GL(V)$ be a   complex reductive group where $\dim V<\infty$, and let $\pi\colon V\to\quot VG$ be the categorical quotient. Let $\NN:=\pi\inv\pi(0)$ be the null cone of $V$, let $H_0$ be the subgroup of $\GL(V)$ which preserves the ideal $\I$ of $\NN$ and let $H$ be a Levi subgroup of $H_0$ containing $G$. We determine the identity component of $H$. In many cases we show that $H=H_0$. For adjoint representations we have $H=H_0$ and we determine $H$ completely. We also investigate the subgroup $G_F$ of $\GL(V)$ preserving a fiber $F$ of   $\pi$ when $V$ is an irreducible cofree $G$-module.  
\end{abstract}

\maketitle
\section{Introduction}
Our base field is $\C$, the field of complex numbers.  Let $V$ be a finite dimensional $G$-module where $G\subset\GL(V)$ is   reductive. Let $R$ denote $\C[V]$. We have the categorical quotient $\pi\colon V\to\quot VG$ dual to the inclusion $R^G\subset R$. Let  $\NN_G:=\pi\inv\pi(0)$ (or just $\NN$) denote  the null cone. Let $G_0=\{g\in\GL(V)\mid f\circ g = f$ for all $f\in R^G\}$. Let $H_0$ denote the subgroup of $\GL(V)$ which preserves $\NN_G$ schematically.  Equivalently, $H_0$ is the group preserving the ideal $\I=R^G_+R$ where $R^G_+$ is the ideal of invariants vanishing at $0$. Let  $G_1$ be a Levi factor of $G_0$ containing $G$ and let $H$ denote a Levi factor of $H_0$ containing $G_1$. We show that $H^0\subset G_1\GL(V)^{G_1}$, hence that $H^0\subset G_1\GL(V)^G$. In many cases  $H_0$ and $G_0$ are reductive, for example, if $V$ is irreducible.
In the case that $V=\lieg$ is a semisimple Lie algebra and $G$ its adjoint group we show that $H=H_0=(\C^*)^r\Aut(\lieg)$ where $r$ is the number of simple ideals in $\lieg$. We also obtain information about the subgroup  of $\GL(\lieg)$ preserving a fiber  of   $\pi$ (other than the zero fiber). We have similar resuts in the case that $V$ is a cofree $G$-module. Our results   generalize those of  Botta,  Pierce and  Watkins \cite{BPW83} and Watkins \cite{Wat82} for the case $\lieg=\lie{sl}_n$. Finally, we show that if $G\subset G'\subset \GL(V)$ where $G'$ is connected reductive such that $\pi$  and $\pi'\colon V\to\quot V{G'}$ have a common fiber, then $R^G=R^{G'}$.

We thank M.\ Ra{\"\i}s for  his help and for the questions and conjectures in his work \cite{Rais1} which led to this paper.

   \section{Equal fibers}
   
 Let $G\subset G'\subset\GL(V)$ be reductive where $G'$ is connected. We have quotient mappings $\pi\colon V\to\quot VG$ and $\pi'\colon V\to\quot V{G'}$. Let $\rho\colon \quot VG\to\quot V{G'}$ denote the canonical map.
  
\begin{theorem} \label{thm:fiber}
 Suppose that there is a fiber $F$ of $\pi$ which is also a fiber of $\pi'$ (as sets). Then $R^G=R^{G'}$.
\end{theorem}

\begin{proof}
 The hypothesis  implies that there is a point $z'\in X':=\quot V{G'}$ such that $\rho\inv(z')$ is a point in $X:=\quot  VG$. Since $\rho$ is surjective,  the minimal dimension of any irreducible component of a fiber is the difference in the dimensions of $X$ and $X'$, so we have that $\dim X=\dim X'$. Then there is a nonempty open subset $U$ of $X'$ such that the fiber  of $\rho$ over any point of $U$ is finite. But for $z'\in X'$, the fiber $(\pi')\inv(z')$ is connected since $G'$ is connected. Hence the fiber  $\rho\inv(z')=\pi((\pi')\inv(z'))$ is connected. It follows that $\rho\colon \rho\inv(U)\to U$ is 1-1 and onto, hence birational. Thus $\rho$ is an isomorphism \cite[II.3.4]{Kr85}
 \end{proof}
 
 \begin{remark} \label{rem:generic} Solomon \cite{Sol05, Sol06} has classified many of the pairs of groups $G\subset G'\subset\GL(V)$ with the same invariants, including the case where $V$ is irreducible. Often, $R^G=R^{G'}$ forces $G=G'$. Suppose that $(V,G)$ is \emph{generic\/}, i.e., it  has trivial principal isotropy groups and the complement of the set of principal orbits has codimension two in $V$. Then $R^G=R^{G'}$ implies that $G=G'$   \cite{Sch07}. 
\end{remark}

\section{Groups preserving the ideal of $\NN$}\label{sec:ideal}

Let $V$ be a $G$-module. We assume that $G$ is a Levi subgroup of $G_0$. Let $H$ be a Levi subgroup of $H_0$ containing $G$. Our aim is to show that $H^0$ is generated by $\GL(V)^G$ and $G^0$.   

\begin{proposition} \label{prop:Gnormal}
Let  $V$, $G$ and $H$ be as above. Then $G$ is normal in $H$.
\end{proposition}

\begin{proof}
Let $p_1,\dots,p_r$ be a set of minimal homogeneous generators of $R^G$.  Let $d_1<d_2<\dots<d_s$ be the distinct degrees of the $p_i$.  Then clearly $H$ preserves the span $W_1$ of the $p_i$ of degree $d_1$. Assuming that  $s>1$, let $W_2'$ be the span of the $p_i$ of degree $d_2$. Then $H$ stabilizes $W_0:=R_{d_2-d_1}W_1$ and $H$ stabilizes $W:=W_2'+W_0=\I\cap R_{d_2}$   where $R_d$ for $d\in\N$ denotes the elements of $R$ homogeneous of degree $d$. Note that $W_2'\cap W_0=W_2'\cap R^G\cdot W_1=0$. Since $H$ is reductive, there is an $H$-stable subspace $W_2$ of $W$ complementary to $W_0$.  Since $G$ acts trivially on $W_2'$, it acts trivially on $W/W_0$ and on $W_2$. Continuing in this way we obtain $H$-modules $W_1,\dots,W_s$ consisting of $G$-invariant functions such that $W':=W_1+\dots + W_s$ generates $R^G$.  Clearly $G$ is the kernel of the action of $H$ on $W'$.
\end{proof}

\begin{corollary}
Suppose that $H_0$ is reductive. Then $G_0$ is reductive and normal in $H_0$.
\end{corollary}

Since $G^0$ is reductive, $H^0$ acts on $G^0$ by inner automorphisms. Hence $H^0=H_1G^0$ where $H_1:=Z_H(G^0)^0$ is the  connected centralizer of $G^0$ in $H$.
\begin{lemma} \label{lemma:theta} Let $g\in G$. Then there is a homomorphism $\theta\colon H_1\to Z(G^0)$ such that
$ghg\inv =\theta(h)h$, $h\in H_1$.
\end{lemma}

\begin{proof}
Let $h\in H_1$. Since conjugation by $h$ preserves the connected components of $G$  there is an element $\theta(h)\in G^0$ such that $hg\inv h\inv=g\inv\theta(h)$.   Let $h_1\in H_1$. Then 
$$
g\inv\theta(h_1h)=h_1hg\inv h\inv h_1\inv=h_1g\inv\theta(h)h_1\inv = h_1g\inv h_1\inv\theta(h)=g\inv \theta(h_1)\theta(h).
$$
Thus $\theta$ is a homomorphism. From $hg\inv h\inv=g\inv \theta(h)$ it follows that $g h g\inv = \theta(h)h$. Since $h$ centralizes $G^0$, so does $g h g\inv$, and we see that $\theta(h)$ centralizes $G^0$. Thus $\theta(h)\in Z(G^0)$.
\end{proof}

\begin{corollary}
Suppose that $G=G_0$ and that $G_0$ is normal in $H_0$. Then $H_0$ is reductive.
\end{corollary}

\begin{proof} As above, we have $(H_0)^0=H_2G^0$ where  $H_2\subset H_0$ is connected and centralizes $G^0$, and $H_0$ is reductive if and only if $H_2$ is reductive. Let $R$ be the unipotent radical of $H_2$. Corresponding to each $g\in G$  there is a homomorphism $\theta\colon H_2\to Z(G^0)$, and since $R$ is unipotent, $\theta(R)=\{e\}$.   Thus $R\subset \GL(V)^G$ where $\GL(V)^G$ is obviously in $H_2$. Thus $R$ is trivial and $H_0$ is reductive.
\end{proof}

Write $H^0=H^0_sG^0_sT$ where $H^0_s$ (resp.\ $G^0_s$) is the semisimple part of $H_1$ (resp.\ $G^0$) and $T:=Z(H^0)^0\subset H_1$ is a torus. Set $T_0:=Z(G^0)^0$.

\begin{corollary}
The group $H^0_s$ is contained in $\GL(V)^G$.
\end{corollary}

\begin{theorem}\label{thm:nullcone} Let $V$, $G$ and $H$ be as above. Then $H^0=\GL(V)^G G^0$. \end{theorem}

\begin{proof} Write $H^0=H^0_sG^0_sT$ as above and set $F:=G/G^0$. Then $F$ normalizes $T$ and by Lemma \ref{lemma:theta}, $F$ acts trivially on $T/T_0$. Thus $T^F$ projects onto $T/T_0$. Choose a  torus $S$ in $(T^F)^0$ complementary to $(T^F\cap T_0)^0$. Then $H^0=H^0_sSG^0$ where $H^0_s S$ lies in $\GL(V)^G$.
\end{proof}

\begin{remark} \label{rem:commutant} Write  $V=\bigoplus_{i=1}^r m_iV_i$ where the $V_i$ are irreducible and pairwise non-isomorphic and $m_iV_i$ denotes the direct sum of $m_i$ copies of $V_i$.   Then the theorem shows that $H^0=G^0\prod_{i=1}^r\GL(m_i)$.
\end{remark}

\begin{example}
Let $\{e\}\neq G\subset\GL(V)$ be finite. Then $\NN_G$, as a set, is just the origin, and it is preserved by  $\GL(V)$. Thus it is essential in Theorem \ref{thm:nullcone} that $H$ preserve $\NN_G$ schematically.
\end{example}

\begin{corollary}\label{cor:nonisomorphic2}
Suppose that $V=\bigoplus_{i=1}^r V_i$ where the $V_i$ are irreducible, nontrivial and pairwise non-isomorphic. Let $H'\subset \GL(V)$ be semisimple. Then the following are equivalent:
\begin{enumerate}
\item $H'\subset H_0$.
\item $H'\subset G_0$.
\end{enumerate}
\end{corollary}

\begin{proposition}\label{prop:irred}
Suppose that $V$ is an irreducible $G$-module. Then $G_0$ and $H_0$ are reductive and $H^0=\C^*G^0$.
\end{proposition}

\begin{proof}
The fixed points of the unipotent radical $R$ of $G_0$ are a $G_0$-stable nonzero subspace of $V$. Thus 
$R$ acts trivially on $V$, i.e., $R=0$. Hence $G_0$ is reductive. Similarly, $H_0$ is reductive.
\end{proof}

\begin{corollary}\label{cor:mW} Suppose that $V=mW$ where $W$ is an irreducible $G$-module. Then $H_0$ is reductive.
\end{corollary}

\begin{proof}
The group $H$ contains $G\times\GL(m)$ which acts irreducibly on $V\simeq W\otimes\C^m$. Thus $H_0$   is reductive.
\end{proof}

In the remainder of this section, we do not assume that $G$ is a Levi subgroup of $G_0$.

\begin{corollary}\label{cor:pWqW^*} Let $G\subset\GL(W)$ and let $V=p W\oplus qW^*$ where $2\leq p\leq q$ and the $G$-modules $W$ and $W^*$ are irreducible  and non isomorphic. Then
\begin{enumerate}
\item $G_0$ and $H_0$ are reductive.
\item $G_0\subset \GL(W)$.
\item $H^0=\GL(p)\GL(q) (G_0)^0$.
\end{enumerate}
\end{corollary}

\begin{proof}
First we consider the case that $G=\GL(W)$. Then   Example \ref{ex:pq} below shows that  $G_0=\GL(W)$ and that $(H_0)^0=\GL(p)\GL(q) \GL(W)$.  
Now the invariants of $\GL(W)$ are generated by those of degree $2$ and  the degree $2$ invariants of    $G$ and of $\GL(W)$ are the same. Thus $G_0$ must be a subgroup of $\GL(W)$ and $(H_0)^0$ must be a subgroup of $\GL(p)\GL(q)\GL(W)$ containing $\GL(p)\GL(q)$. Hence $(H_0)^0=\GL(p)\GL(q)H_1$ where $H_1\subset\GL(W)$. Note that $GH_1$ is a finite extension of $H_1$. Since $W$ is an irreducible $G$-module and $G_0$ and $GH_1$ contain $G$, both $G_0$ and $H_1$ (hence $(H_0)^0)$  are reductive and we have (1) and (2).   Theorem \ref{thm:nullcone} gives (3).
\end{proof}

\begin{lemma}\label{lem:H0reductive}
Suppose that $V^G=(0)$ and let $V=\bigoplus_{i=1}^r m_i V_i$ be the isotypic decomposition of $V$ where the $V_i$ are pairwise non-isomorphic $G$-modules.  Suppose that $\lieh_0(m_iV_i)\subset m_iV_i$ for all $i$. Then $H_0$ is reductive.
\end{lemma}

\begin{proof} For any $i$,   $G(H_0)^0$ is a finite extension of $(H_0)^0$ which contains $G\prod_i\GL(m_i)$. The latter group acts irreducibly on $m_iV_i$, hence the image of $G(H_0)^0$ in $\GL(m_i V_i)$ is reductive for all $i$.  It follows that $(H_0)^0$ is reductive, hence that $H_0$ is reductive.
\end{proof}

 \begin{corollary}\label{cor:nonisomorphic}  
Suppose that $V_i$ is an irreducible   nontrivial  $G_i$-module where $G_i$ is reductive and $\C[V_i]^{G_i}\neq\C$, $i=1,\dots,r$.     Let $V:=\bigoplus_i m_iV_i$ with the canonical action of $G:=G_1\times\dots\times G_r$ where $m_i\geq 1$ for all $i$. Then  $H_0$ is reductive.
\end{corollary} 

\begin{proof}
Suppose that $\lieh_0$ is not contained in $\bigoplus_i\End(m_i V_i)$. Since  $\lieh_0$ is $H$-stable, it must contain one of the irreducible $G_i\times\GL(m_i)\times G_j\times\GL(m_j)$-modules   $\Hom(m_iV_i,m_jV_j)$, $i\neq j$. Without loss of generality suppose that $\lieh_0\supset \Hom(m_2V_2,m_1V_1)$.
Let $f\in\O(m_1V_1)^{G_1}$ be a  nonconstant  homogeneous invariant of minimal degree $d\geq 2$.
Let $\phi\in\Hom(m_2V_2,m_1V_1)$. Then $\phi$ sends $f$ to the function $h(v_1,v_2):=df(v_1)(\phi(v_2))$ where $v_i\in m_iV_i$, $i=1$, $2$. Clearly there is  a $\phi$ such that $h\neq 0$.  Thus $h$ is a nonzero element of bidegree $(d-1,1)$ in $\C[m_1V_1\oplus m_2V_2]$. But by the minimality of $d$ and the fact that no nonzero  invariant in $\C[m_2V_2]$ has degree 1, there is no element of $\I$ of this bidegree. Hence $\Hom(m_2V_2,m_1V_1)$ does not preserve $\I$, a contradiction. 
Thus $\lieh_0$ is contained in $\bigoplus_i\End(m_iV_i)$ and one can apply Lemma \ref{lem:H0reductive}.
\end{proof}

\begin{corollary}
Suppose that $G\subset\GL(V)$ is a finite group generated by pseudoreflections. Then $H_0$ is reductive.
\end{corollary}
\begin{proof}
We have that $V=\bigoplus V_i$ and $G=\prod G_i$ where $G_i\subset\GL(V_i)$ is an irreducible group generated by pseudoreflections. Now apply Corollary \ref{cor:nonisomorphic}.
\end{proof}

\begin{proposition}\label{prop:orthogonal}  Suppose that $V$ is an orthogonal representation of $G$ where $V^G=(0)$. Then $H_0$ is reductive.
 \end{proposition} 

\begin{proof}  We have an isotypic decomposition $V=\bigoplus_i m_i V_i\bigoplus n_j (W_j\oplus W_j^*)$ where the $V_i$ are irreducible nontrivial orthogonal representations of $G$ and the $W_j$ are irreducible nonorthogonal representations of $G$.  Note that for each $i$ there is a quadratic invariant $p_i\in\C[m_iV_i]^G$ and for each $j$ a quadratic invariant (a contraction) $q_j\in\C[n_j(W_j\oplus W_j^*)]^G$. Suppose that $\lieh_0$ is not contained in $\bigoplus_i\End(m_iV_i)\bigoplus_j\End(n_j(W_j\oplus W_j^*))$. For example, suppose that there is a nonzero element $\phi$ of $\lieh_0$ whose restriction to $m_2V_2$ has nonzero projection to $m_1V_1$. Then we have the function $h(v_1,v_2):=dp_1(v_1)(\phi(v_2))$ for $v_1\in m_1V_1$ and $v_2\in m_2V_2$. As before, the actions of $G$ and the $\GL(m_i)$ guarantee that we can assume that $h\neq 0$. Now the bidegree of $h$ is $(1,1)$ and $h\in\I$. However, there are no nonconstant invariants of bidegree $(a,b)$ in $\C[m_1V_1\oplus m_2V_2]$ for $a\leq 1$ and $b\leq 1$. Thus $h$ cannot lie in $\I$. One similarly gets contradictions for all the possible ways that  $\lieh_0\not\subset  \bigoplus_i\End(m_iV_i)\bigoplus_j\End(n_j(W_j\oplus W_j^*))$ can occur. Finally, note that the normalizer $N$ of the image of $G$  in $\GL(n_j(W_j\oplus W_j^*))$ contains an element interchanging the copies of $W_j$ and $W_j^*$. Thus $N$  acts irreducibly and we can now apply the argument of Lemma \ref{lem:H0reductive}.
 \end{proof}

\begin{corollary}
If $G$ is any one of the following groups, then $H_0$ is reductive for any representation $V$ of $G$ with $V^G=(0)$.
\begin{enumerate}
\item $\SO(n)$, $n\geq 3$.
\item $\mathsf{G_2}$, $\mathsf{F_4}$, $\mathsf{E_8}$.
\item $\mathsf{B}_{4n+3}$ and $\mathsf{B}_{4n+4}$, $n\geq 0$.
\item $\mathsf{D}_{4n}$, $n\geq 1$.
\end{enumerate}
\end{corollary}

\section{Some examples and a conjecture}

We give examples where $G_0$ is not reductive and we give examples where $G_0$ is reductive but $H_0$ is not.

\begin{example}
Let $V$ be a $G$-module and $W$ a $G$-module such that $\O(V\oplus W)^G=\O(V)^G$. Then $\Hom(V,W)$ is contained in the radical of  $\lieg_0$  so that $G_0$ and $H_0$ are not reductive. A concrete example is given by $G=\SL_4$ and $V\oplus W=\wedge^2\C^4\oplus\C^4$ with the obvious $G$ action.
\end{example}

\begin{example}
Let $W$ be an irreducible $G$-module where  $W^G=(0)$ and $\O(W)^G\neq \C$. Let $V=W\oplus\C$ where $G$ acts trivially on $\C$. Then $\lieg_0\subset \lie{gl}(W)$ while
$\Hom(\C,W)$ is contained in the Lie algebra of the radical of $H_0$.
\end{example}

\begin{example}\label{ex:pq}
Let $1\leq p\leq q$ and consider the $G=\GL(W)$ representation  on $V=pW\oplus qW^*$ where $W=\C^n$, $n\geq 1$. (See Corollary \ref{cor:pWqW^*}.)\ By classical invariant theory, the $G$-invariants are just the contractions of elements of the copies of $W$ with elements of the copies of $W^*$. Let $U$ denote $W\oplus W^*\simeq\C^{2n}$.

Three cases arise:

\smallskip
\noindent Case 1:  $p=q=1$. Then our invariant is the bilinear form $(\ ,\ )$ corresponding to the matrix $J:=\sm 0 I I 0\subset\GL(2n)$, i.e., $(x,y)=x^tJy$, $x$, $y\in U$.  Thus $G_0=\Orth(2n)$ and $H_0=\C^*G_0$.
 
\smallskip
\noindent Case 2: $p=1$, $q>1$. Then $H_0$ contains a copy of $\GL(q)$ and the action of $H_0$ on the invariants is a representation $H_0\to\GL(q)$ whose kernel is $G_0$. Thus $H_0=\GL(q)G_0$. A matrix computation shows that $G_0=\GL(W)\ltimes(\wedge^2(W^*)\otimes\C^q)$. If $x\in W$ and $y_1,\dots, y_q\in W^*$, then the unipotent radical of $G_0$ sends $(x,y_1,\dots,y_q)$ to $(x,y_1+B_1x,\dots,  y_q+B_q x)$  where for each $j$, $B_j$ is a skew symmetric matrix, $B_j\in\wedge^2(W^*)\subset \Hom(W,W^*)$.

\smallskip
\noindent Case 3: $p\geq 2$.  We show that $G_0=\GL(W)$, that $H_0=H$ and that $H^0=\GL(p)\GL(q)\GL(W)$ and we determine $H$. First suppose that $p=q=2$. Then $G_0$ preserves the inner products on $2U$, i.e., $G_0$ is a subgroup of $\Orth(2n)$. Moreover, $G_0$ preserves the skew product on $2U$ sending $x$, $y$ to $x^tKy$ where $K=\sm 0 I {-I}0$. Hence $G_0$ lies in the intersection of $\Orth(2n)$ and $\Sp(2n)$ which is the copy of $\GL(W)$ acting on $U$ by the matrices $\sm A 0 0 {{}^t A\inv}$, $A\in\GL(W)$.  Clearly, as long as $2\leq p\leq q$ we must have that $G_0=G=\GL(W)$. We have a representation $\phi\colon H_0\to\GL(pq)$ given by the action of $H_0$ on the $pq$ generators of the invariants. The kernel of $\phi$ is   $G_0=G$. Thus $H_0$ is reductive. By Theorem \ref{thm:nullcone} we have $H^0=\GL(p)\GL(q)\GL(W)$. Let $h\in H$.   If $h$   stabilizes $pW$ and $qW^*$, then $h$ induces an automorphism of $\GL(W)$ which is trivial on $\C^*I$ and must be   inner on $\SL(W)$. Hence modulo an element of $\SL(W)$, $h$ lies in the centralizer of $\GL(W)$, which is $\GL(p)\GL(q)$. Hence $h\in H^0$. The only other possibility is that $h$ interchanges the copies of $pW$ and $qW^*$. This can only happen if $p=q$. Thus $H$ is connected if $p\neq q$ and $H/H^0$ has order two if $p=q$.
\end{example}

\begin{example}\label{ex:finite}
Let $G=\Z/4\Z\subset\C^*$ and let $V=\C^2$ where $\xi(a,b)=(\xi^2a,\xi b)$ for $(a,b)\in\C^2$, $\xi\in G$. Since $G$ is finite, $G_0=G$. Let $x$ and $y$ be the usual coordinate functions on $V$. Then the invariants are generated by $x^2$, $xy^2$ and $y^4$. Consider the  element $\phi\in\End(V)$ which sends $(a,b)$ to $(0,a)$ for $a$, $b\in\C$. Then $\phi$ acts on $\C[V]$ by the derivation $x\partial/\partial y$. This derivation preserves $\I$ and it follows that $\phi$ is a basis of the Lie algebra of the unipotent radical of $H_0$.
\end{example}

\begin{example} \label{ex:cstar}
Let $G=\C^*$ and let $V$ be the $p+q+r$ dimensional representation with weights $-1$ of multiplicity $p$, $1$ of multiplicity $q$ and $2$ of multiplicity $r$ where $p$, $q$, $r\in\N$ and $pqr\neq 0$. If $x_i$, $y_j$ and $z_k$ are corresponding coordinate functions, then the invariants are generated by the monomials $x_iy_j$ and $x_ix_{i'}z_k$. We have $G_0=G$ while the radical of $H_0$ has Lie algebra spanned by  the linear mappings corresponding to the derivations $y_j\partial/\partial z_k$.
\end{example}

\begin{example}\label{ex:sl2}
Let $V\oplus W=S^2(\C^n)\oplus \C^n$ with the obvious action of $G=\SL_n$, $n\geq 2$. Then using classical invariant theory \cite{Sch78a} one computes that the invariants have homogeneous generators  $p$ and $q$ of bidegrees $(n,0)$ and $(n-1,2)$, respectively. Now $\Hom(V,W)$ contains a copy of $W^*$ where $\xi\in W^*$ sends $v\in V$ to $i_\xi(v)\in W$ (contraction). Then this copy of $W^*$ acts on $\C[V\oplus W]$ sending a polynomial $f(v,w)$ into $df(v,w)(0,i_\xi(v))$, $v\in V$, $w\in W$. This action annihilates $p$ and sends $q$ to a subspace of $\O(V\oplus W)$ of bidegree $(n,1)$ transforming under $G$ as $W^*$. But the only way to get a copy of $W^*$ in this bidegree is to multiply $p$ times the copy of $W^*$ in degree 1 in $\O(V\oplus W)$. Thus $\I$ is preserved. It is now easy to establish that the unipotent radical of $H_0$ has Lie algebra the copy of $W^*$ in $\Hom(V,W)$.  
\end{example}

\begin{conjecture}
 If $G$ is semisimple and $V$ is generic (see \ref{rem:generic}) with $V^G=(0)$, then $H_0$ is reductive.
\end{conjecture}

\section{Cofree Representations}

   Recall that $V$ is   \emph{cofree\/} if $R$ is a free module over $R^G$. Equivalently, $R^G$ is a polynomial ring and $\pi \colon V\to\quot  VG$   is equidimensional \cite[17.29]{Sch80}. If $p_1,\dots,p_d$ are minimal homogeneous generators of $R^G$, then we can identify $\pi$ with the polynomial map $p=(p_1,\dots,p_d)\colon V\to \C^d$.  Cofreeness is equivalent to the fact that the $p_i$ form a regular sequence in $\C[V]$.
  See \cite{Sch78b} for the classification of cofree representations of the simple algebraic groups and \cite{Lit89} for the classification of irreducible cofree representations of semisimple algebraic groups.
   
   We say that $ G'\subset\GL(V)$   stabilizes a fiber $F$ of $\pi$ if $G'$ preserves $F$  schematically, i.e., preserves the ideal $I_F$ of $F$.   
     
     \begin{proposition} \label{prop:fibercofree} Suppose that   $G$  is reductive and $V$ is a cofree   $G$-module. 
If $G'\subset\GL(V)$ stabilizes a fiber of $\pi\colon V\to\quot VG$, then
$G'$ stabilizes $\NN_G$.
     \end{proposition}
 
 \begin{proof} Let $F$ be a fiber of $\pi$. Then there are constants $c_i$, $i=1,\dots,d$, such that $I_F$ is the ideal generated by $p_i-c_i$, $i=1,\dots,d$. Let $0\neq f\in I_F$ and let $\gr f$ denote the nonzero homogeneous part of $f$ of largest degree. Then the elements $\gr f$ for $0\neq f\in I_F$ generate a homogeneous ideal $I$ which obviously contains $\I$. We show that $I\subset \I$ so that $I=\I$. If $G'$ preserves $I_F$, it preserves $I=\I$, and we have the proposition.
 
 Let $d_i$ be the degree of $p_i$, $i=1,\dots,d$. Let $0\neq f\in I_F$ where $\gr f$ is homogeneous of degree $r$. We have $f=\sum a_i(p_i-c_i)$ where $a_1,\dots,a_n\in R$. Let $s=\max_i\{\deg a_i+d_i\}$. Let $a_i'$ denote the homogeneous part of $a_i$ of degree $s-d_i$. If $s>r$, then we must have that $\sum_i a_i'p_i=0$. Since the $p_i$ are a regular sequence, this relation is generated by the Koszul relations $p_jp_i-p_ip_j=0$, $1\leq i<j\leq d$. Hence there are $b_{ij}\in R$, $b_{ij}=-b_{ji}$, such that
$$
\sum_i a_i'(p_i-c_i)= \sum_{i\neq j} b_{ij} (p_j(p_i-c_i)-p_i(p_j-c_j)) =\sum_{i\neq j} b_{ij} (c_j(p_i-c_i)-c_i(p_j-c_j)) 
$$ 
where for fixed $i$, $\deg\sum_{j\neq i} b_{ij}c_j<s-d_i$. Thus we may replace each $a_i$ by a polynomial of degree less than $s-d_i$ without changing $f$.  Continuing inductively we  reduce to the situation that $\deg a_i\leq r-d_i$ for all $i$. Let $a_i'$ denote the homogeneous degree $r-d_i$ term in $a_i$, $i=1,\dots,d$. Then $\gr f = \sum_i a_i'p_i \in\I$.
 \end{proof}

\begin{example}
Let $G=\C^*$ and $V=\C^3$ with coordinate functions  $x$, $y$ and $z$ corresponding to weights $-1$, $1$ and $2$.    The fiber defined by xy=1 and $x^2z=0$ is the fiber  defined by  $xy=1$ and $z=0$, and it has a symmetry which interchanges $x$ and $y$. However, this is not a symmetry of the ideal generated by the invariants. Thus Proposition \ref{prop:fibercofree} does not hold in case the representation is not cofree.
\end{example}

\begin{remark} \label{rem:smooth}
Let $F$ be a principal fiber of $\pi$ where $V$ is cofree. Then $d\pi$ has rank $d=\dim\quot VG$ on $F$ so that $F$ is smooth. It follows that $G'$ preserves $I_F$ if and only if $G'$ preserves the set $F$. 
\end{remark}
 
 \begin{corollary} \label{cor:three} Let $V=\bigoplus_{i=1}^r V_i$ where the $V_i$ are pairwise non-isomorphic nontrivial $G$-modules and $V$ is cofree.
 Suppose that $G\subset G'\subset\GL(V)$ where $G'$ is connected semisimple.  Then the following are equivalent.
  \begin{enumerate}
  \item $R^G=R^{G'}$.
  \item $G'$ preserves a fiber of $\pi\colon V\to\quot VG$.
\item $G'$ preserves $\NN_G$.
\end{enumerate}
 \end{corollary}
 
 \begin{proof} Use Corollary \ref{cor:nonisomorphic2} and  Proposition \ref{prop:fibercofree}.
  \end{proof}

\begin{corollary} \label{cor:cofree}
Let   $V$ be an irreducible nontrivial cofree $G$-module with $R^G\neq \C$. Let $F\neq\NN$ be a fiber of $\pi\colon V\to\quot VG$ and let $G_F$ be the subgroup of $\GL(V)$ stabilizing $F$.  Then
\begin{enumerate}
\item $G_0\subset G_F\subset H_0$ are reductive.
\item $H^0=\C^*(G_0)^0$.
\item $G_F/G_0$ is finite.
\end{enumerate}
\end{corollary}

\begin{proof}
Parts (1) and (2) are clear. Since $F\neq\NN$, it is only stabilized by  a finite subgroup of $\C^*$, hence we have (3).
\end{proof}

 It would be nice to find an example of an irreducible module $V$ of a semisimple group $G$ with $G=(G_0)^0$  such that the subgroup of $\GL(V)$ fixing a fiber $F$ of $\pi\colon V\to\quot VG$, $F\neq\NN_G$, has dimension bigger than $\dim G$.  

\begin{remark} Let $V$ be an irreducible nontrivial cofree representation of a simple algebraic group $G$ such that $R^G\neq\C$.
The cases for which $G\neq (G_0)^0$ are as follows (we use the numbering and notation of \cite{Sch78a}).
\begin{enumerate}
\item $(\phi_3,\mathsf{B_3})$.
\item $(\phi_4,\mathsf{B_4})$.
\item $(\phi_5,\mathsf{B_5})$.
\item $(\phi_1,\mathsf{G_2})$.
\end{enumerate}
\end{remark}

  \section{The adjoint case}   
  
  Let $\lieg$ be a simple Lie algebra.
 Choose a Cartan subalgebra $\liet$ of $\lieg$ and a base $\Pi$ of the root system. Choose $x_\alpha\in \lieg_\alpha$ and $y_\alpha\in\lieg_{-\alpha}$, $\alpha\in\Pi$, such that $(x_\alpha$, $y_\alpha$, $[x_\alpha,y_\alpha])$ is a standard $\lie{sl}_2$ triple.
Then there is a unique order 2 automorphism $\psi$ of $\lieg$ which is $-1$ on $\liet$ and sends $x_\alpha$ to $-y_\alpha$, $\alpha\in\Pi$. 
  
 Now let $\lieg=\lieg_1\oplus\dots\oplus\lieg_r$ where the $\lieg_i$ are simple ideals of the Lie algebra $\lieg$. Let $\psi_i\in\Aut(\lieg_i)$ be as above. Let $G$ denote the adjoint group of $\lieg$ and let $G_0$, $H_0$ and $H$ be as in the introduction.
 
\begin{theorem} \label{thm:main1} We have that $H=(\C^*)^r\Aut(\lieg)$ and that $G_0\simeq(\Z/2\Z)^rG$ where the $i$th copy of $\Z/2\Z$ is generated by $-\psi_i$.
\end{theorem}

 \begin{proof} By Dixmier \cite{Dix79}, $(G_0)^0=G$, and using Corollary \ref{cor:nonisomorphic} we obtain that $H_0=H$ where $H^0=(\C^*)^r G$. Hence if $\phi\in H$, we obtain an automorphism $\sigma$ of $\lieg\simeq\ad\lieg\subset\lieh\subset\lie{gl}(\lieg)$ where $\ad(\sigma(X))=\phi\circ\ad X\circ\phi\inv$, $X\in\lieg$.  Clearly $\Aut(\lieg)\subset H$, so replacing $\phi$ by $\phi\circ\sigma\inv$ we can arrange that $\phi\circ\ad X\circ\phi\inv=\ad X$ for all $X\in\lieg$. Then by Schur's lemma, $\phi\in(\C^*)^r\subset H^0$ so that $H=(\C^*)^r\Aut(\lieg)$. If we start with $\phi\in G_0$, then since $\phi$ induces the identity on $\C[\lieg]^G$, so does $\sigma$, and it follows from Schur's lemma that $\phi$ is a product $\prod_i\lambda_i \sigma_i$ where,  for all $i$,  $\sigma_i\colon \lieg_i\to\lieg_i$ is an automorphism and $\lambda_i\in\C^*$ acts via multiplication on $\lieg_i$. But $\phi$ has to preserve the invariants of degree 2 of each $\lieg_i$, hence $\lambda_i=\pm 1$ for all $i$. Now   \cite[Theorem 2.5]{Sch07} shows that,  for each $i$, $\lambda_i\sigma_i\in G_i$ or $\lambda_i\sigma_i \in(-\psi_i) G_i\neq G_i$, where $G_i$ is the adjoint group of $\lieg_i$. Hence $G_0\simeq(\Z/2\Z)^rG$.
  \end{proof}

\begin{corollary}[\cite{BPW83}]
Let $\lieg=\lie{sl}_n$. Then $H$ is generated by $G$,  $\C^*$ and transposition.
\end{corollary}

\begin{proof}
In the case of $\lie{sl}_n$  with the usual choice of $\liet$ and $\Pi$, the automorphism $\psi$ is $X\mapsto -X^t$, $X\in\lie{sl}_n$.  Then $\psi$ generates the group of outer automorphisms of $\lie{sl}_n$ (which is the trivial group for $n=2$). Hence $H$ is generated by $G$, $\C^*$ and transposition.
\end{proof}

\begin{corollary}[\cite{Wat82}]
Let $G_F$ be the subgroup of $\GL(\lie{gl}_n)$ which preserves the $G:=\PGL(n)$-orbit $F$ of an element    $x_0$ of  $\lie{gl}_n$ which has nonzero trace and distinct eigenvalues. Then $G_F$ is generated by $G$ and transposition.
\end{corollary}
\begin{proof} The condition on $x_0$ implies that $F$ is a smooth  fiber of the quotient mapping (see Remark \ref{rem:smooth}).
 Write $x_0=\mu I+y_0$ where $\mu\in\C^*$, $y_0\in\lie{sl}_n$ and $I$ is the $n\times n$ identity matrix. Then $F$ is just $\mu I+F_1$ where $F_1=G\cdot y_0$. We may write an element of $G_F$ as $\sm 1 0 c{ \lambda g}$ where $c\in\lie{sl}_n$, $g\in\GL(\lie{sl}_n)$ is a linear mapping preserving the schematic null cone of $\lie{sl}_n$ and $\lambda\in\C^*$  (use \ref{prop:irred}  and \ref{prop:fibercofree}). Then $g$ is in $G$ or $g$ is an element of $G$ composed with transposition. Applying the inverse of $g$ we obtain an element $h$ of the form $y\mapsto \lambda y+c$, $y\in F_1$. We need to show that $c=0$. Suppose not. Let $g\in G$ such that $gc\neq c$. Then $h\inv ghg\inv(y)=  y+c'$, $0\neq c'\in\lie{sl}_n$,  $y \in F_1$. Thus $F_1=F_1+c'$. It follows that for any   invariant polynomial $p$ on $\lie{sl}_n$, $p(y+nc')=p(y)$ for all $y\in F_1$ and $n\in\Z$. Thus $dp(y)(c')=0$ for any $y\in F_1$. But the covectors $dp(y)$ for $y\in F_1$ span $(\lie{sl}_n)^*$. Thus $c'=0$, a contradiction.
 \end{proof}


\newcommand{\noopsort}[1]{} \newcommand{\printfirst}[2]{#1}
  \newcommand{\singleletter}[1]{#1} \newcommand{\switchargs}[2]{#2#1}
  \def\cprime{$'$}
\providecommand{\bysame}{\leavevmode\hbox to3em{\hrulefill}\thinspace}
\providecommand{\MR}{\relax\ifhmode\unskip\space\fi MR }
\providecommand{\MRhref}[2]{%
  \href{http://www.ams.org/mathscinet-getitem?mr=#1}{#2}
}\providecommand{\href}[2]{#2}

\end{document}